# PARAMETRIC BOOTSTRAP APPROXIMATION TO THE DISTRIBUTION OF EBLUP AND RELATED PREDICTION INTERVALS IN LINEAR MIXED MODELS

By Snigdhansu Chatterjee, Partha Lahiri and Huilin Li

*University of Minnesota, University of Maryland
and University of Maryland*

Empirical best linear unbiased prediction (EBLUP) method uses a linear mixed model in combining information from different sources of information. This method is particularly useful in small area problems. The variability of an EBLUP is traditionally measured by the mean squared prediction error (MSPE), and interval estimates are generally constructed using estimates of the MSPE. Such methods have shortcomings like under-coverage or over-coverage, excessive length and lack of interpretability. We propose a parametric bootstrap approach to estimate the entire distribution of a suitably centered and scaled EBLUP. The bootstrap histogram is highly accurate, and differs from the true EBLUP distribution by only $O(d^3 n^{-3/2})$, where $d$ is the number of parameters and $n$ the number of observations. This result is used to obtain highly accurate prediction intervals. Simulation results demonstrate the superiority of this method over existing techniques of constructing prediction intervals in linear mixed models.

**1. Introduction.** Large scale sample surveys are usually designed to produce reliable estimates of various characteristics of interest for large geographic areas. However, for effective planning of health, social and other services, and for apportioning government funds, there is a growing demand to produce similar estimates for smaller geographic areas and for other subpopulations. To meet this demand, it is necessary to supplement the survey data with other relevant information that is often obtained from different administrative and census records. In many small area applications, mixed linear models are now routinely used in combining information from various sources and explaining different sources of errors. These models incorporate area specific random effects which explain the "between small area variations," not otherwise explained by the fixed effects part of the model.

For a good review on small area and linear mixed model research, the readers are referred to the book by Rao (2003), and two recent review papers by Rao (2005) and Jiang and Lahiri (2006). Several other applications of linear mixed models may be found in McCulloch and Searle (2001). Point prediction using the empirical best linear unbiased predictor (EBLUP) and the associated mean square prediction error (MSPE) estimation have been studied extensively. See Jiang, Lahiri and Wan (2002), Rao (2005) and Jiang and Lahiri (2006) for a review on the subject, especially on the latest development on resampling methods for MSPE estimation. However, little progress has been made outside the basic study of the first two moments, for example, on the properties of quantiles (central or tail) of predictors, or on the effect of high dimensionality of the parameters.

For example, research on interval estimates in small area studies is typically limited to some special cases of the Fay–Herriot model (described in detail in Section 2), where the traditional estimates are of the form $EBLUP \pm z_{\alpha/2}\sqrt{mspe}$. Here $mspe$ is an estimate of the true $MSPE$ of the EBLUP, and $z_{\alpha/2}$ is the upper $100(1-\alpha/2)\%$ point of the standard normal distribution. The coverage probabilities of such intervals may converge to the nominal level $1-\alpha$; but the intervals are not efficient, in the sense they have either under-coverage or over-coverage problem, depending on the particular choice of the MSPE estimator. More precisely, the coverage error of such interval is of the order $O(n^{-1})$ or higher, which is not accurate enough for most applications of small area studies, many of which involve small sample size $n$.

In this paper we address the problem of approximating the distribution of a predictor, and applying it to obtain prediction intervals, in a very general framework of linear mixed models. We consider the following model from Das, Jiang and Rao (2004):

$$(1.1) \quad \mathbf{Y}_n = \mathbf{X}\beta + \mathbf{Z}\mathbf{v}_q + \mathbf{e}_n,$$

where $\mathbf{Y}_n \in \mathbb{R}^n$ is a vector of observed responses, $\mathbf{X}_{n\times p}$ and $\mathbf{Z}_{n\times q}$ are known matrices and $\mathbf{v}_q$ and $\mathbf{e}_n$ are independent random variables with dispersion matrices $D_q(\psi)$ and $R_n(\psi)$, respectively. Here $\beta \in \mathbb{R}^p$ and $\psi \in \mathbb{R}^k$ are fixed parameters.

The mixed ANOVA model, and the longitudinal models including the Fay–Herriot model and the nested error regression model are special cases of (1.1). We can consider both balanced and unbalanced lay-outs in the above framework. In addition, we develop our theory and methodology allowing for the parameter dimension $d = p + k$ to grow with sample size $n$. Dimension dependent asymptotics are extremely important in the current context, since many small areas may have sample sizes comparable to dimensions of the regression and variance components parameters; see, for example, Jiang (1996) for their use and importance in linear mixed models.

Our approach toward approximating the distribution of a predictor is to employ parametric bootstrap. We concentrate on the *empirical best linear predictor* (EBLUP), owing to its wide popularity and use. We establish in Theorem 3.1 that the bootstrap histogram incurs an error of $O(d^3 n^{-3/2})$ in approximating the distribution of a centered and scaled EBLUP. For estimating the distribution of centered and scaled *estimators*, under standard regularity conditions and fixed $d$, the normal approximation based on the central limit theorem has an error of $O(n^{-1/2})$, but the bootstrap can achieve higher-order accuracy, with typical approximation error of $O(n^{-1})$. Theorem 3.1 may be seen as an extension of this higher-order accuracy phenomenon, in the context of prediction. Although our motivation and terminology comes from small area context, our bootstrap methodology and theoretical results are directly applicable to other usages of mixed linear models.

There are several potential applications of a highly accurate approximation of the entire distribution of the EBLUP. For example, it may be used to obtain (a) bagging predictors, (b) computing mean squared errors or other risks, (c) hypothesis testing, (d) calibration of traditional estimators, and (e) prediction interval construction. In this paper, we concentrate on the last application, since prediction intervals combine features of both point prediction and hypothesis testing nicely, and have not been extensively explored in small area or other mixed linear model contexts.

Prediction intervals are useful in small area studies in several ways. For example, prediction intervals may help establish if different counties have similar resources and needs, or if different ethnic or other subpopulation groups are equally exposed to a particular disease. Our simultaneous concentration on dimension asymptotics is also relevant. It has long been recognized that health, economic activity and other measures of human well-being depend on a number of exogenous and endogenous factors, many of which must be measured at the individual level and incorporated in the model. In statistical terms, this translates to high dimensionality of $\beta$ and $\psi$.

In Section 2, we review some of the existing techniques for predictor distribution approximation and interval estimate construction. We pay special attention to the usage of resampling in such approximations/constructions. For prediction intervals, available literature is heavily concentrated on special cases of the Fay–Herriot model. Since traditional intervals perform poorly in terms of coverage or length or both, many attempts have been made to fine tune and calibrate them, often using resampling. To the best of our knowledge, approximation of the entire distribution of a predictor has not been attempted in general small area problems, and we briefly review the related research for independent data.

In Section 3, we present our bootstrap algorithm for the Das, Jiang and Rao model (1.1). Our main result is that the sup-norm distance between the distribution of EBLUP and its bootstrap approximation is $O(d^3 n^{-3/2})$. A direct corollary is that the bootstrap prediction interval has coverage accuracy of $O(d^3 n^{-3/2})$. Note that our proposed prediction interval is a *bootstrap interval*, which is different from the traditional approaches of obtaining asymptotic intervals first and then calibrating it. Our interval can be calibrated one or more times to achieve coverage accuracy of $O(d^5 n^{-5/2})$ or higher, if needed.

We performed several simulation experiments in order to study how our percentile bootstrap interval estimate compares with existing techniques. A sample of these studies are reported in Section 4. The main message from the simulations is that the prediction intervals resulting from the proposed parametric bootstrap perform considerably better than the traditional techniques, which is a reflection of the high order accuracy theoretically established in Section 3.

## 2. A review of predictor distribution approximation and interval construction.

2.1. *Approximating distributions of predictors.*  Considerable theoretical research has been carried out in the prediction of a random variable that is independent of $\mathbf{Y}_n$, and has density $\xi(\cdot|\beta, \psi)$. In terms of expected Kullback–Leibler divergence, the naive plug-in predictor density $\xi(\cdot|\hat{\beta}, \hat{\psi})$ performs poorly compared to Bayesian predictors $\int \xi(\cdot|\beta, \psi) \pi(\beta, \psi|\mathbf{Y}_n)$, see for example, Aitchison (1975), Murray (1977), Ng (1980), Komaki (1996, 2001, 2006) and George, Liang and Xu (2006). Harris (1989) showed that the *bootstrap predictor*

$$(2.1) \qquad \xi^*(\cdot) = \int \xi(\cdot|s, t) \, d\mathcal{L}^*(s, t)$$

also performs better than the naive plug-in predictor. Recently, Fushiki, Komaki and Aihara (2004) have shown that the bootstrap predictor (2.1) is asymptotically equivalent to a Bayesian predictor with Hartigan's $M$-prior. The $M$-prior has certain optimality properties which may be found in Hartigan (1964, 1998). In a related work, Fushiki, Komaki, Aihara (2005) show that the Harris predictor is related to bagging of Breiman (1996).

In small area or other mixed linear model contexts, the random variable of interest depends on $\mathbf{Y}_n$, unlike the framework described above. Also, performance measures other than expected Kullback–Leibler divergence may be of interest, for example, length and coverage of prediction intervals.

2.2. *A review of interval estimation techniques.* For a general mixed linear model, Jeske and Harville (1988) proposed a prediction interval for a mixed effect, but did not include the effect of estimated unknown variance components on the accuracy of their proposed interval.

Jiang and Zhang (2002) used a distribution-free method for constructing prediction intervals for a future observation under a non-Gaussian linear mixed model, based on the theory developed by Jiang (1998). This technique does not employ any area specific information and can be useful in constructing intervals when there is no survey data on the response variable. Jiang and Zhang (2002) proposed another method which can be applied to the situation when the sample size is large within each area. This is a technique of first obtaining the EBLUP for the random effects and the residuals. Then, under conditions sufficient to imply that the number of times each random effect is repeated (i.e., number of observations in each small area) tends to infinity, the empirical distribution of random effects as well as the residuals converge appropriately. This technique fails when we do not have large samples for each small area, a situation that is common in many small area applications.

Recently, Hall and Maiti (2006b) have studied parametric bootstrap for general mixed models in several aspects, including interval estimation. A review of their approach toward interval estimation may be found in Rao (2005). In Section 3, we discuss in detail how their model, results and asymptotics differ from ours.

Other than the above three papers, research on small area prediction intervals is largely concentrated on special cases of the Fay–Herriot model, described below:

1. Conditional on $\theta = (\theta_1, \ldots, \theta_n)^T$, $\mathbf{Y}_n = (Y_1, \ldots, Y_n)^T$ follows a $n$-variate normal distribution with mean $\theta$ and dispersion matrix $D$ with *known* diagonal entries $D_i > 0$ and off-diagonal entries 0. Here (and in the sequel) all vectors are taken to be column vectors, for any vector (matrix) $a$ ($A$), the notation $a^T$ ($A^T$) denotes its transpose.
2. The variable $\theta$ follows a $n$-variate normal distribution with mean $\mathbf{X}\beta$ for a known $n \times p$ matrix $\mathbf{X}$ and unknown but fixed vector $\beta \in \mathbb{R}^p$. The dispersion matrix is $A\mathbf{I}_n$, where the matrix $\mathbf{I}_n$ is the $n$ dimensional identity matrix and $A$ is an unknown constant.

There are several options for constructing interval estimates for $\theta_i = \mathbf{x}_i^T \beta + v_i$. One may use only the Level 1 model for the observed data, or only the Level 2 model for the borrowed strength component, or a combination of both. The interval for $\theta_i$ based only on the Level 1 model is given by $I_i^D(\alpha): Y_i \pm z_{\alpha/2} D_i^{1/2}$, where $z_{\alpha/2}$ is the $(1-\alpha/2)$th standard normal quantile. Obviously, for this interval, the coverage probability is $1 - \alpha$. However, it is not efficient, since its average length is too large to make any reasonable conclusion. This is due to the high variability of the point predictor $Y_i$.

An interval based only on Level 2 ignores the crucial area specific data that is modeled in Level 1, and hence falls short on two counts: it fails to be relevant to the specific small area under consideration, and it fails to achieve sufficient coverage accuracy. A small example given later in this section demonstrates this latter property.

Thus, interval estimation techniques that combine both levels of the Fay–Herriot model are required. A popular approach is to employ empirical Bayes methodology. Cox (1975) proposed the following empirical Bayes interval:

$$I_i^C(\alpha) : (1 - \hat{B}_i)Y_i + \hat{B}_i \mathbf{x}_i^T \hat{\beta} \pm z_{\alpha/2} D_i^{1/2}(1 - \hat{B}_i)^{1/2},$$

where $\hat{B}_i$ and $\hat{\beta}$ are estimators of $B_i = D_i/(A + D_i)$ and $\beta$, respectively, and $\mathbf{x}_i^T$ is the $i$th row of $\mathbf{X}$. Under standard regularity conditions, $\mathbb{P}(\theta_i \in I_i^C(\alpha)) = 1 - \alpha + O(n^{-1})$, where $\mathbb{P}$ denotes a probability measure induced by the joint distribution of Level 1 and Level 2. Thus, this prediction interval attains the desired coverage probability asymptotically, but the coverage error is of order $O(n^{-1})$, which is not accurate enough for many small area applications. This lack of accuracy may partially be due to the additional variability resulting from estimation of $\beta$ and $A$. Currently, MSPE estimators are available in several mixed linear models, see for example, Jiang Lahiri and Wan (2002), Datta, Rao and Smith (2005), Hall and Maiti (2006a). Naive empirical Bayes intervals constructed using EBLUP, MSPE estimators and standard normal quantiles typically have an error of $O(n^{-1})$ or higher.

For a special case of the Fay–Herriot model with common mean and equal sampling variances $D_i = D$, Morris (1983a) incorporated the additional uncertainty due to the estimation of the hyperparameters. However, Basu, Ghosh and Mukerjee (2003) showed that the resulting empirical Bayes interval proposed by Morris (1983a) still has coverage error of $O(n^{-1})$. They used analytical calibration of the Morris' interval to reduce the coverage error to $o(n^{-1})$. They also showed that with suitable analytical approximations in place, an interval due to Carlin and Louis (1996), page 98, and a new interval, have coverage error of the order $o(n^{-1})$. Datta et al. (2002) used similar analytical calibration in a more general Fay–Herriot model, and obtained a prediction interval with coverage error of $O(n^{-3/2})$. Morris (1983b) considered a variation of his (1983a) work with the use of a hierarchical Bayes type point estimator. Hill (1990) suggested a general framework which, in the Fay-Herriot setting matches with an exact hierarchical Bayes confidence interval. Datta et al. (2002) followed up Hill's idea to obtain an interval with coverage error of $O(n^{-1})$.

Apart from the analytical approaches, calibration using different bootstrap techniques has been popular. The methods differ in the generation of the bootstrap samples and the type of correction made. For a special case of the Fay–Herriot model where $Y_1, \ldots, Y_n$ are independent and identically distributed, Laird and Louis (1987) proposed three different resampling strategies: (a) usual nonparametric bootstrap by sampling with replacement from the data, (b) a semi-parametric method, assuming density at the first level of their two level model is known but that at the second level is unknown, and (c) the parametric bootstrap. In mixed linear models, the nonparametric and semi-parametric bootstrap approximation of the distribution of the EBLUP are generally not consistent. Once the bootstrap sample (nonparametric, semi-parametric or parametric) is generated, the next challenge is to find a method that corrects the empirical Bayes confidence intervals $I_i^C(\alpha)$ to achieve better coverage. Laird and Louis (1987) considered an imitation of the hierarchical Bayes approach.

Carlin and Gelfand (1990, 1991) point out that the hierarchical Bayesian methods like those of Laird and Louis (1987) lead to a *lengthening* of the empirical Bayes interval, which is not the same as a *correction*. They discuss an example where increasing the length further exacerbates the coverage bias. They suggest parametric bootstrap to calibrate the empirical Bayes interval.

Calibration of intervals has been one of the major uses of bootstrap for some time, and can lead to considerable improvement of coverage accuracy. Coupled with use of bias correction, use of pivotal or nearly pivotal statistics, and Edgeworth corrections, improvements from calibration can sometimes be dramatic. See Abramovitch and Singh (1985), Beran (1990a, 1990b), the book by Efron and Tibshirani (1993) and references therein for further details on these issues. On the other hand, calibration is both time and computational effort intensive, often requiring iterative searches; it typically increases variability; the results often lack straightforward interpretability; and successive calibrating steps typically have diminishing returns in terms of improvement of coverage. It is not always clear what property of an interval, that is, length, coverage, end points or some other characteristic, ought to be calibrated, see for example, DiCiccio and Efron (1996) and the discussions of it by Hall and Martin (1996), Lee and Young (1996); and the interesting example in the rejoinder. Some calibrating options do not exist for multivariate confidence or prediction regions. Asymptotic results suggest calibrated intervals have better coverage accuracy, but do not consider the variability induced by the calibration, do not represent performance in finite samples; or reflect the degree in which the finite sample results depend on unknown parameters and their estimators. Nevertheless, calibration is an excellent tool to improve coverage of intervals; though it seems sensible to use a more accurate interval and little or no calibration; rather than a less accurate interval with intensive calibration. The bootstrap interval we obtain in Section 3 is one such highly accurate interval, and requires the same amount of computational effort as one round of bootstrap based calibration of Carlin and Gelfand (1990, 1991).

Hall (2006) suggested an application of the nonparametric bootstrap confidence interval based on the generated $\theta_i^\star$'s only. In the small area context, this may be applicable when the differences between the small areas are minor, or carried only in the fixed effects. In surveys, robustness is always an important issue, and the practitioners are always interested in efficient nonparametric methods. However, due to scarce data at the small area level, nonparametric estimators tend to under-perform, often severely. This is because the nonparametric models typically permit the generation of bootstrap histograms based on a synthetic model or the regression model, but do not permit approximation of the conditional distribution of $\theta_i$ given the data $\mathbf{Y}_n$. As a result, the nonparametric bootstrap prediction interval for $\theta_i$ is likely to underweight the area specific data. Accurate weighting of the area specific data is important for achieving good coverage properties, as the example below shows. Hall (2006) also pointed out the importance of parametric bootstrap in small area estimation and other related problems.

EXAMPLE. Consider the following special case of the Fay–Herriot model where $\sigma_i \equiv 1$, and $\mathbf{x}_i^T \beta \equiv \mu$. Thus, at Level 1, $Y_i$'s given the $\theta_i$'s are independently distributed as $N(\theta_i, 1)$ random variables; and at Level 2, the $\theta_i$'s are independent, identically distributed as $N(\mu, \tau^2)$ random variables. The estimators of $\mu$ and $\tau^2$ are given, respectively, by $\hat{\mu} = \bar{Y}$, $\hat{\tau}^2 = \max(0, s^2 - 1)$, where $s^2 = \sum(y_i - \bar{y})^2/(n-1)$. Assume $\hat{\tau}^2 > 0$, a condition that is satisfied in many problems. The bootstrap procedure would require us to generate $\theta_i^* \stackrel{\text{iid}}{\sim} N[\hat{\mu}, \hat{\tau}^2]$ and $Y_i^* | \theta_i^* \stackrel{\text{ind}}{\sim} N[\theta_i^*, 1]$. Then we have $\hat{\mu}^* = \bar{Y}^*$, $\hat{\tau}^{2*} = \max(0, s^{*2} - 1)$ where $s^{*2} = \sum(y_i^* - \bar{y}^*)^2/(n-1)$. An obvious Level 2 based bootstrap prediction interval for $\theta_i$ that is *not* area specific, is given by

$$(2.2) \qquad (\hat{\mu} - t_1\sqrt{\hat{\tau}^2}, \hat{\mu} + t_2\sqrt{\hat{\tau}^2}),$$

where $(t_1, t_2)$ are cutoff points satisfying $\mathbb{P}(\hat{\mu}^* - t_1\sqrt{\hat{\tau}^{2*}} \leq \theta^* \leq \hat{\mu}^* + t_2\sqrt{\hat{\tau}^{2*}}) = 1 - \alpha$.

It can be shown that interval (2.2) has coverage of $1 - \alpha + O(n^{-1/2})$ which makes it consistent, but hardly accurate enough. The lack of accuracy is due to the use of the Level 2 distribution only, so that the Level 1 data $Y_i$ plays no special role in the interval construction.

In Bayesian terminology, the Level 2 of the Fay–Herriot model essentially corresponds to a *prior* on $\theta_i$, while the Level 1 model yields the *likelihood*. Using only the "prior knowledge" (Level 2 distribution) does not even yield consistency in general. However, in some instances using the Level 2 distribution in conjunction with bootstrap can have a calibration effect that obtains $O(n^{-1/2})$ consistency, as shown above.

**3. Parametric bootstrap prediction interval for a general linear mixed model.** We consider the model:

$$\mathbf{Y}_n = \mathbf{X}\beta + \mathbf{Z}\mathbf{v}_q + \mathbf{e}_n, \tag{3.1}$$

where $\mathbf{X}$ is a known $(n \times p)$ matrix, $\mathbf{Z}$ is a known $(n \times q)$ matrix, $\mathbf{Y}_n \in \mathbb{R}^n$ is the vector of observed data, $\beta \in \mathbb{R}^p$ is a fixed but unknown parameter vector, and $\mathbf{v}_q \in \mathbb{R}^q$ and $\mathbf{e}_n \in \mathbb{R}^n$ are random variables following the normal distributions $N_q(0, D_q)$ and $N_n(0, R_n)$, respectively. The integer $q$ may depend on $n$, thus $q \equiv q_n$. Assume the sequence $\{\mathbf{v}_q\}$ and $\{\mathbf{e}_n\}$ are independent. The first term $\mathbf{X}\beta$ represents the fixed effects, and the second term $\mathbf{Z}\mathbf{v}_q$ the random effects. Thus $\mathbf{X}\beta + \mathbf{Z}\mathbf{v}_q$ constitute the *signal* component of the observed data, while $\mathbf{e}_n$ is the *noise*. The properties of the signal are of interest, which depend on the unknown parameters $\beta$, $D_q$ and $R_n$.

Assume that the $(q \times q)$ matrix $D_q$ and the $(n \times n)$ matrix $R_n$ are known up to a $(k \times 1)$ vector of unknown parameters, thus $D_q = D_q(\psi)$ and $R_n = R_n(\psi)$ for a fixed but unknown $\psi = (\psi_1, \ldots, \psi_k)^T \in \mathbb{R}^k$. Note that the dispersion matrix of the observed data $\mathbf{Y}_n$ is given by

$$\Sigma_n = \Sigma_n(\psi) = R_n(\psi) + \mathbf{Z}D_q(\psi)\mathbf{Z}^T.$$

We henceforth drop the $n$ from $\mathbf{Y}_n$, $\mathbf{e}_n$, $R_n$ and $\Sigma_n$, and $q$ from $\mathbf{v}_q$ and $D_q$ to simplify notation. We take $d = p + k$, the dimension of the parameter space. Let $\theta = (\beta, \psi)$ denote the unknown parameters.

Das, Jiang and Rao (2004) show that several linear mixed models, including analysis of variance (ANOVA) models and longitudinal models of both balanced and unbalanced nature are special cases of the model (3.1). Unbalanced ANOVA models arise, for example, when $R = \sigma_0^2 \mathbf{I}_n$; and $D = \text{diag}(\sigma_1^2 \mathbf{I}_{r_1}, \ldots, \sigma_{k-1}^2 \mathbf{I}_{r_{k-1}})$ where $\mathbf{I}_r$ is the $r \times r$ identity matrix. Here $\psi$ is the vector of variance components $\psi = (\sigma_0^2, \ldots, \sigma_{k-1}^2)$. Unbalanced longitudinal models arise when $\Sigma$ has a block diagonal structure.

Let $T = c^T(\mathbf{X}\beta + \mathbf{Z}\mathbf{v})$, where $c$ is a fixed and known $(n \times 1)$ vector. The case where $c$ is a $n \times m$ matrix obtains multidimensional predictive quantities, and their treatment is similar to the univariate case described below, with some minor algebraic variations. We concentrate on univariate $T$ for easier exposition. The conditional distribution of $T$ given $\mathbf{Y}$ is $N(\mu_T, \sigma_T^2)$, where

$$\begin{aligned}\mu_T &= c^T\mathbf{X}\beta + c^T\mathbf{Z}D\mathbf{Z}^T\Sigma^{-1}(\mathbf{Y} - \mathbf{X}\beta) \\ &= c^T R\Sigma^{-1}\mathbf{X}\beta + c^T\mathbf{Z}D\mathbf{Z}^T\Sigma^{-1}\mathbf{Y}\end{aligned} \tag{3.2}$$

and

$$\sigma_T^2 = c^T\mathbf{Z}(D - D\mathbf{Z}^T\Sigma^{-1}\mathbf{Z}D)\mathbf{Z}^T c. \tag{3.3}$$

Generally, $\beta$ and $\psi$ (and hence $D$ and $R$) are estimated from the data $\mathbf{Y}$ by using the marginal distribution of $\mathbf{Y}$, given by $N_n(\mathbf{X}\beta, \Sigma)$. The resulting estimates $\hat{\mu}_T$ and $\hat{\sigma}_T$ of the mean and variance of $T$ are expressions similar to (3.2) and (3.3), with $\hat{\beta}$ and $\hat{\psi}$ in place of $\beta$ and $\psi$.

For algebraic simplicity, in the rest of this paper we assume that $\mathbf{X}$ is full column rank and use the estimator $\hat{\beta} = (\mathbf{X}^T\mathbf{X})^{-1}\mathbf{X}^T\mathbf{Y}$. This is the *ordinary least squares* estimator of $\beta$. Using other estimators like the maximum likelihood estimator or the *weighted least squares* estimator, with appropriate conditions on the weights, is another possibility. This makes little difference in the asymptotic analysis as long as the weights are smooth functions of $\psi$. Estimator $\hat{\psi}$ of $\psi$ is typically obtained by maximum likelihood or restricted maximum likelihood techniques.

Based on the fact that $\sigma_T^{-1}(T - \mu_T)$ is a standard normal pivot, the traditional approach to interval estimation for $T$, reviewed in Section 2, is to take $(\hat{\mu}_T \pm z\sqrt{mspe})$ for some estimator *mspe* of *MSPE* and the appropriate Normal quantile $z$. Unfortunately, $\hat{\sigma}_T^{-1}(T - \hat{\mu}_T)$ is not a pivot, and the traditional approach produces too short or too long intervals. Let the distribution of $\hat{\sigma}_T^{-1}(T - \hat{\mu}_T)$ be $\mathcal{L}_n$. Recognizing that $\mathcal{L}_n$ is not the standard normal distribution, we propose to estimate it using parametric bootstrap.

Define
$$\mathbf{Y}^* = \mathbf{X}\hat{\beta} + \mathbf{Z}\mathbf{v}^* + \mathbf{e}^*$$
where $\mathbf{v}^* \sim N_q(0, D(\hat{\psi}))$ and $\mathbf{e}^* \sim N_n(0, R(\hat{\psi}))$ are independent of each other.

From $\mathbf{Y}^*$, obtain $\hat{\beta}^*$ and $\hat{\psi}^*$ using the same techniques used to obtain $\hat{\beta}$ and $\hat{\psi}$ earlier. Next, obtain $\hat{\mu}_T^*$ and $\hat{\sigma}_T^*$ using $\hat{\beta}^*$ and $\hat{\psi}^*$ using (3.2) and (3.3). Define $T^* = c^T(\mathbf{X}\hat{\beta} + \mathbf{Z}\mathbf{v}^*)$. The distribution of
$$\hat{\sigma}_T^{-1*}(T^* - \hat{\mu}_T^*),$$
conditional on the data $\mathbf{Y}$, is the parametric bootstrap approximation $\mathcal{L}_n^*$ of $\mathcal{L}_n$. Using this approximation, we then proceed to obtain the interval estimate for $T$ as $(\hat{\mu}_T + q_1\hat{\sigma}_T, \hat{\mu}_T + q_2\hat{\sigma}_T)$, where $q_1$ and $q_2$ are appropriate quantiles of the bootstrap approximation $\mathcal{L}_n^*$ of $\mathcal{L}_n$.

Our main result is that $\mathcal{L}_n^*$ approximates $\mathcal{L}_n$ up to $O(d^3 n^{-3/2})$ terms. In order to state the assumptions for our result, let us introduce some terminology and notation now. For any function $f(\psi): \mathbb{R}^a \to \mathbb{R}$, $f'(\psi)$ denotes its first derivative written as a $a \times 1$ column vector; $f''(\psi)$ denotes the $a \times a$ second derivative matrix. For a symmetric matrix A, $\lambda_{\max}$ and $\lambda_{\min}$, respectively, denote its maximum and minimum eigenvalue.

The following are the assumptions for our result in this section:

1. The following relations hold:

(3.4) $$\|\mathbf{X}^T c\| = O(1),$$

$$\|\mathbf{X}^T \Sigma^{-1} \mathbf{Z} D \mathbf{Z}^T c\| = O(1), \tag{3.5}$$

$$c^T \mathbf{Z} D \mathbf{Z}^T c = O(1), \tag{3.6}$$

$$c^T \mathbf{Z} D \mathbf{Z}^T \Sigma^{-1} \mathbf{Z} D \mathbf{Z}^T c = O(1). \tag{3.7}$$

In addition,

$$\sigma_T^2 = c^T \mathbf{Z}(D - D\mathbf{Z}^T \Sigma^{-1} \mathbf{Z} D)\mathbf{Z}^T c > M > 0,$$

for some constant $M > 0$.

2. Assume that

$$\sup_{1 \leq i \leq n} \sum_{j=1}^{p} \left[ \sum_{a=1}^{n} X_{ja} \Sigma_{ai}^{1/2} \right]^2 = O(p/n), \tag{3.8}$$

$$\lambda_{\min}(n^{-1} \mathbf{X}^T \mathbf{X}) > M > 0, \tag{3.9}$$

for some constant $M > 0$.

3. The eigenvalues of the matrices $D$ and $R$ lie in $(L^{-1}, L)$ for some $L > 1$. The eigenvalues of $D(\hat{\psi})$ and $R(\hat{\psi})$ lie in $(L^{-1}/2, 2L)$. The eigenvalues of $\Sigma$ lie in a compact set on the positive half of the real line.

In the representations

$$D = D_1(\psi) D_1^T(\psi), \qquad \hat{D} = D_1(\psi) \Lambda_D(\hat{\psi}) D_1^T(\psi), \tag{3.10}$$

$$R = R_1(\psi) R_1^T(\psi), \qquad \hat{R} = R_1(\psi) \Lambda_R(\hat{\psi}) R_1^T(\psi), \tag{3.11}$$

where $\Lambda_R$ and $\Lambda_D$ are diagonal matrices, the following conditions are satisfied:

All the entries of the $q \times q$ matrix $\Lambda_D = \operatorname{diag}(\Lambda_{D1}, \ldots, \Lambda_{Dq})$ and the $n \times n$ matrix $\Lambda_R = \operatorname{diag}(\Lambda_{R1}, \ldots, \Lambda_{Rn})$ have three bounded continuous derivatives.

We denote by $\Lambda'_D$ the $k \times q$ matrix whose $(j, i)$th entries are given by

$$((\Lambda'_D))_{j,i}(\psi) = \frac{\partial}{\partial \psi_j} \Lambda_{Di}(\psi), \qquad j = 1, \ldots, k; \ i = 1, \ldots, q.$$

The $(j, i)$th entry of the $k^2 \times q$ matrix $\Lambda''_D$ is

$$((\Lambda''_D))_{j,i}(\psi) = \frac{\partial^2}{\partial \psi_{j_1} \partial \psi_{j_2}} \Lambda_{Di}(\psi), \qquad j_1 + (j_2 - 1)k = j,$$

$$j_1, j_2 = 1, \ldots, k, \ j = 1, \ldots, k^2, \ i = 1, \ldots, q.$$

The $(j, i)$th entry of the $k^3 \times q$ matrix $\Lambda_D^{(3)}$ is

$$((\Lambda_D^{(3)}))_{j,i}(\psi) = \frac{\partial^3}{\partial \psi_{j_1} \partial \psi_{j_2} \partial \psi_{j_3}} \Lambda_{Di}(\psi),$$

where $j_1 + (j_2 - 1)k + (j_3 - 1)k^2 = j$, $j_1, j_2, j_3 = 1, \ldots, k$, $j = 1, \ldots, k^3$, and $i = 1, \ldots, q$.

We define the $k \times n$ matrix $\Lambda'_R$, the $k^2 \times n$ matrix $\Lambda''_R$ and the $k^3 \times n$ matrix $\Lambda_R^{(3)}$ along identical lines as above.

The following conditions are assumed:

$$\text{(3.12)} \quad \lambda_{\max} \Lambda_D'^T(\psi) \Lambda_D'(\psi) = O(1),$$

$$\text{(3.13)} \quad \lambda_{\max} \Lambda_R'^T(\psi) \Lambda_R'(\psi) = O(1),$$

$$\text{(3.14)} \quad \lambda_{\max} \Lambda_D''^T(\psi) \Lambda_D''(\psi) = O(1),$$

$$\text{(3.15)} \quad \lambda_{\max} \Lambda_R''^T(\psi) \Lambda_R''(\psi) = O(1),$$

$$\text{(3.16)} \quad \lambda_{\max} \Lambda_D^{(3)T}(\psi^*) \Lambda_D^{(3)}(\psi^*) < M = O(1),$$

$$\text{(3.17)} \quad \lambda_{\max} \Lambda_R^{(3)T}(\psi^*) \Lambda_R^{(3)}(\psi^*) < M = O(1),$$

for some constant $M > 0$ for all $\psi^*$ in a neighborhood of the true value $\psi$.

4. Let $S = (k/n)^{1/2}(\hat{\psi} - \psi)$. Assume that all the moments of $\|S\|$ are $O(1)$. Moreover, the following relations are also satisfied:

$$\text{(3.18)} \quad \mathbb{E}S_j = O(\sqrt{k/n}), \quad j = 1, \ldots, k,$$

$$\text{(3.19)} \quad \mathbb{E}S_a S_b = O(\sqrt{k/n}), \quad a, b = 1, \ldots, k,$$

$$\text{(3.20)} \quad \mathbb{E}S_j (\mathbf{Z}\mathbf{v} + \mathbf{e})_i = O(\sqrt{k/n}), \quad j = 1, \ldots, k, \; i = 1, \ldots, n,$$

$$\text{(3.21)} \quad \mathbb{E}S_a S_b (\mathbf{Z}\mathbf{v} + \mathbf{e})_i = O(\sqrt{k/n}), \quad a, b = 1, \ldots, k, \; i = 1, \ldots, n.$$

We now state our main theorem for this section.

THEOREM 3.1. *Under the Assumptions* (1)–(4), *if* $d^2/n \to 0$, *we have*

$$\text{(3.22)} \quad \sup_{q \in \mathbb{R}} |\mathcal{L}_n(q) - \mathcal{L}_n^*(q)| = O_P(d^3 n^{-3/2}).$$

The proof of Theorem 3.1 is given in the Appendix. A direct application of Theorem 3.1 is the following result on highly accurate prediction intervals.

THEOREM 3.2. *Under the Assumptions* (1)–(4) *and* $d^2/n \to 0$, *for any* $\alpha \in (0, 1)$, *if* $q_1$ *and* $q_2$ *are real numbers such that*

$$\mathcal{L}_n^*(q_2) - \mathcal{L}_n^*(q_1) = 1 - \alpha,$$

*we have*

$$\text{(3.23)} \quad \mathbb{P}[\hat{\mu}_T + q_1 \hat{\sigma}_T \leq T \leq \hat{\mu}_T + q_2 \hat{\sigma}_T] = 1 - \alpha + O(d^3 n^{-3/2}).$$

Theorem 3.2 follows directly from Theorem 3.1, hence we omit its proof. Since the Fay–Herriot model (described in Section 2) is an important example, we state the results for it in a separate corollary below.

COROLLARY 3.1. *In the Fay–Herriot model, assume that the matrix* $\mathbf{X}$ *is full column rank, the diagonal entries* $h_{ii}$ *of the projection matrix on the columns of* $\mathbf{X}$ *satisfy* $\sup_i h_{ii} = O(p/n)$, *the Level 1 variances* $\{D_i\}$ *lie in a compact subset of* $(0, \infty)$, *and the estimator* $\hat{A}$ *of* $A$ *is positive. Then, for any* $i \in \{1, \ldots, n\}$, *if* $\hat{\theta}_i^{EB} = (1 - \hat{B}_i) Y_i + \hat{B}_i \mathbf{x}_i^T \hat{\beta}$, $\hat{\theta}_i^{EB*} = (1 - \hat{B}_i^*) Y_i^* + \hat{B}_i^* \mathbf{x}_i^T \hat{\beta}^*$, *we have*

$$\begin{aligned}(3.24)\quad &\mathbb{P}[\theta_i \in (\hat{\theta}_i^{EB} + q_{i1} D_i^{1/2} (1 - \hat{B}_i)^{1/2}, \hat{\theta}_i^{EB} + q_{i2} D_i^{1/2} (1 - \hat{B}_i)^{1/2})] \\ &= 1 - \alpha + O(p^3 n^{-3/2});\end{aligned}$$

*where* $\hat{B}_i = D_i / (\hat{A} + D_i)$, *and* $(q_{i1}, q_{i2})$ *satisfy*

$$\begin{aligned}&\mathbb{P}^*[\theta_i^* \in (\hat{\theta}_i^{EB*} + q_{i1} D_i^{1/2} (1 - \hat{B}_i^*)^{1/2}, \hat{\theta}_i^{EB*} + q_{i2} D_i^{1/2} (1 - \hat{B}_i^*)^{1/2})] \\ &= 1 - \alpha + O_P(p^3 n^{-3/2}).\end{aligned}$$

The notation used in Corollary 3.1 are standard ones, that is, $\mathbb{P}^*$ is the probability on the resampling scheme conditional on the data, $\hat{B}_i^* = D_i / (\hat{A}^* + D_i)$, where $\hat{\beta}^*$ and $\hat{A}^*$ are the estimators computed on the bootstrap data $\mathbf{Y}^*$. Here conditional on the data, $\theta_i^* \sim N(\mathbf{x}_i^T \hat{\beta}, \hat{A})$, and $Y_i^* | \theta_i^* \sim N(\theta_i^*, D_i)$ independently. Corollary 3.1 is easily derived from Theorem 3.2, and we omit the details of its proof. A slightly different approach to the same result may be found in the unpublished manuscript Chatterjee and Lahiri (2002). We now discuss the assumptions leading to our main result Theorem 3.1, and some additional features of our result.

REMARK 1 (On the dimension of the random effect vector). Note that the dimension $q$ of the random effect $\mathbf{v}$ is arbitrary which may or may not depend on $n$. Owing to this generalization, our analysis is for $T = c^T(\mathbf{X}\beta + \mathbf{Z}\mathbf{v})$, rather than the more traditional $\tilde{T} = c_1^T \beta + c_2^T \mathbf{v}$. Since $\mathbf{X}$ is full column rank, the fixed effects in $T$ and $\tilde{T}$ are equivalent.

REMARK 2 (On the technical assumptions). In the development of all the assumptions above, we have preferred simplicity over generality. The requirement $d^2 n^{-1} \to 0$ is standard in dimension asymptotics. Assumption 1 is in order to ensure $T$ as a nontrivial quantity, that is, it ensures that both the fixed component and the variance of the random component of $\mu_T$ are $O(1)$, and the variance $\sigma_T^2$ is bounded away from zero and infinity. By suitably scaling the norm of the vector $c$ this assumption is satisfied.

Assumption 2 is a standard assumption on the behavior **X**. It ensures that the norm of each fixed effects covariate is of suitable order, and the fixed effects design is not singular. This assumption can be modified to suit cases where **X** is not full column rank, but such generalizations are routine.

Assumption 3 is on standard differentiability and eigenvalue conditions. Here again, we have tried to adopt simple conditions rather than the most general ones. Note that the existence of the representations (3.10) and (3.11) are not part of the assumptions, and these representations will be established in the proof of Theorem 3.1.

Also note that the eigenvalues of $D(\hat{\psi})$ and $R(\hat{\psi})$ are estimates of the variance components in typical applications. Note that we do not allow these to be zero, since these must always lie in $(L^{-1}/2, 2L)$. However, $L$ may be arbitrarily large, consequently this assumption does not limit the applicability of our results.

In Assumption 4 we take all moments of $S$ to exist in order to achieve simplicity. Our result involves computation of several terms involving $S$, and having all the moments of $S$ available simplifies the algebra. In most applications, both $\psi$ and $\hat{\psi}$ lie in a compact set, hence this is not a strong condition. The other moment conditions on $S$ given by (3.18)–(3.19) are routine. These hold when $\hat{\psi}$ is obtained using either maximum likelihood or restricted maximum likelihood formulation, see Jiang (1998) for related developments.

Conditions (3.20)–(3.21) are interesting, since they effectively set a limit to the amount of dependency structure we can have in $\Sigma$. In order to visualize this, suppose $\hat{\psi}^{(-i)}$ is the estimator of $\psi$ obtained by using only those observations that are independent of $Y_i$; and let $S^{(-i)} = (k/n)^{1/2}(\hat{\psi}^{(-i)} - \psi)$. Then, a sufficient condition for $\mathbb{E}S_j(\mathbf{Zv} + \mathbf{e})_i = O((k/n)^{1/2})$ is that $S - S^{(-i)} = O_P((k/n)^{1/2})$.

This is routinely achieved, and in particular, if $Y_i$ is independent of all but a finite number of observations, we have $S - S^{(-i)} = O_P((k/n)^{1/2})$. This is the typical situation is almost all applications of small area studies. Thus, the effect of Assumption 4 is to restrict the complexity of the matrices $D$ and $R$.

REMARK 3 (On the nature of prediction intervals). A prominent application of the highly accurate approximation of $\mathcal{L}_n(\cdot)$ by $\mathcal{L}_n^*(\cdot)$ is stated in Theorem 3.2, that is, in the construction of prediction intervals. Note that these are *bootstrap intervals*, as opposed to the traditional intervals described in Section 2, some of which are improved with bootstrap corrections.

However, Theorem 3.2 does not describe the nature of the bootstrap prediction intervals, since the choice of $q_1$ and $q_2$ can be quite arbitrary.

These may be chosen to ensure either an *equal tail* property of the interval; whereby $\mathcal{L}_n^*(q_1) = \alpha/2$ and $\mathcal{L}_n^*(q_1) = 1 - \alpha/2$; or we may chose these according to a *minimum length of interval* property, that is, we minimize the length $\hat{\sigma}_T(q_2 - q_1)$. The simulation experiments reported in Section 4 show that both equal tailed and minimum length bootstrap prediction intervals typically achieve the desired coverage accuracy without the use of elaborate calibrations; and the minimum length interval is always slightly shorter than the equal tailed one.

REMARK 4 (On multivariate prediction). Note that in place of the real valued $T$ studied above, we could have a vector valued $T$ with little change in methodology. The algorithmic and algebraic details are similar, and the main result of high order accuracy of distributional approximation (3.22) holds. The major difference between univariate and multivariate prediction is in the construction of prediction regions. Instead of the two points $q_1$ and $q_2$, we need to obtain probability concentration regions from the bootstrap distribution. Such regions can be obtained using various data depth notions and shape features, for example, as in Yeh and Singh (1997). This is a separate issue from the one addressed in this paper, and will be handled in a different paper. Note that multidimensional probability concentration regions can be quite hard to calibrate in practice. Some techniques, such as calibration of the end points of an interval, are not available in this case.

REMARK 5 (Asymptotics on total sample size $n$). One important feature of Theorem 3.1 is that the asymptotic limits are obtained with total sample size $n$ tending to infinity. The total sample size $n$ is the sum total of all observations made, counting each repeated measurement on each individual unit in each small area as a distinct observation. This allows Theorem 3.1 to be used with considerable flexibility, for example, when number of individual units in small areas are large, or when number of small areas are large, or both. However, requirements of asymptotic negligibility, as in (3.20)–(3.21), must still be met. Our assumptions are designed for the more realistic applications where number of small areas are large.

In general, for mixed linear models asymptotic limits are obtained either when the number of observations in each small area tends to infinity, or when the number of small areas tend to infinity; see McCulloch and Searle (2001) and Rao (2003) for details. Theorem 3.1 is a breakthrough, owing to the greater flexibility it allows in asymptotics.

REMARK 6 (On area specific properties). The *area-specific signal* for each small area is of $T_i = \mathbf{x}_i^T \beta + \mathbf{Z}_i^T \mathbf{v}$, conditional on the observed $i$th small area data $Y_i$. Distributions of predictors for such area specific signals are effectively captured by our bootstrap predictive distribution approximation.

Consequently, our bootstrap prediction intervals for $T_i$ are also *area specific*. Extensions to compare two or more small areas can be obtained by similar techniques, see comment on multidimensional prediction above.

In the prediction interval described in Theorem 3.2, we have considered unconditional coverage, where probabilities are computed over the joint distribution of $\mathbf{Y}$ and $\mathbf{v}$. This establishes the performance of the area-specific interval that depends on $\mathbf{v}$ conditional on $\mathbf{Y}$, as well as variability due to observations $\mathbf{Y}$.

Alternatively, one might compute the *area-specific* (*random*) *coverage*, which is defined as $\mathbb{P}[T_i \in I_{Pi}|\mathbf{Y}_i]$, where $I_{Pi}$ is the prediction interval. The interval proposed in Theorem 3.2 achieves $O_P(d^2/n)$ order of area-specific coverage accuracy, since some smoothing effects arising from the distribution of $\mathbf{Y}$ are absent. This is no worse (and in some cases, better) than the area-specific coverage obtained by other techniques in special cases of the general linear model (3.1).

REMARK 7 (On calibration). Both the unconditional as well as the area-specific coverage can be improved by calibration. The use of calibration coupled with resampling is an active topic of research, and some discussion on this has been presented in Section 2. The coverage accuracy of the prediction interval of Theorem 3.2 can be improved to $O(d^5 n^{-5/2})$ with one round of calibration, and further still with more calibration. Such calibration may be done either on the probabilities corresponding to the two end points as in DiCiccio and Efron (1996), or on the true coverage of the interval. Some of our simulations, not reported in this paper, suggest that it is not always beneficial to attempt boosting the theoretical coverage probability, disregarding other properties of the interval. For example, variability of calibrated intervals are greater than uncalibrated ones, minimum length property is almost never preserved, and the results are quite dependent on the parameters and fixed constants of the problem. Hence, it seems reasonable to work with a good predictive distribution as in Theorem 3.1, instead of starting with a naive interval and embarking on intense iterative calibration.

REMARK 8 [The parallel work of Hall and Maiti (2006b)]. Recently, Hall and Maiti (2006b) studied parametric bootstrap methods for general small area models, and considering the overlap of the topics studied in their paper and this one, deserve special mention. For this comment, we use some notation from Hall and Maiti (2006b) whenever they are not in conflict with the notations in the rest of this paper, but use our notation otherwise.

For a suitable function $f_i(\beta)$ involving co-variates $X_i = (X_{i1}, \ldots, X_{in_i})$ and parameter $\beta$, they consider random effect $\Theta_i \sim Q(\cdot; f_i(\beta), \xi)$, and conditional on $\Theta_i$, the data $Y_{ij}$ are independent observations from $R(\cdot; \psi(\Theta_i), \eta_i)$, for $j = 1, \ldots, n_i$, $i = 1, \ldots, m$. Here $\psi(\cdot)$ is a known link function, $\xi$ and $\eta_i$'s are either parameters or known constants, and $Q(\cdot)$ and $R(\cdot)$ are known probability distribution functions. They go on to study calibration of the mean squared prediction error (MSPE) and interval estimation with parametric bootstrap.

Their model is broad enough to handle nonlinear mixed effects, which our model (1.1) cannot do. However, their assumption of $Y_{ij}$'s being independent means that they do not consider longitudinal models, or other models with temporal or spatial dependence. This is essentially the case $R$ being a multiple of the identity matrix in our set-up. Our model is broader than Hall and Maiti's in including several varieties of dependence structure.

The interval estimate from their Section 2.8 is

$$\hat{I}_\alpha = \mathbf{x}_i^T \hat{\beta} \pm z_{\alpha/2} \hat{A}^{1/2} \tag{3.25}$$

for the Fay–Herriot model. Rao (2005) noted that this interval does not make use of the area-specific direct estimator, unlike the prediction interval proposed by Chatterjee and Lahiri (2002). Hall and Maiti (2006b) calibrate this interval for better coverage accuracy, improving from their result

$$\mathbb{P}[\Theta_i \in \hat{I}_\alpha] = 1 - \alpha + O(m^{-1}). \tag{3.26}$$

The result (3.26) hold when the probability statement is on the marginal distribution of random effect $\Theta_i$, and estimators $\hat{\beta}$ and $\hat{A}$ are independent of $i$th area data $Y_{i1}, \ldots, Y_{in_i}$.

Our probability statements in Theorems 3.1 and 3.2 are, however, on the *joint variability* of the random effects and data $(\Theta_i, Y_{i1}, \ldots, Y_{in_i})$. Also note that Theorem 3.2 is obtained as $n = \sum_{i=1}^m n_i \to \infty$, while (3.26) is obtained as $m \to \infty$. Since some of the $n_i$'s can be large, the speed of convergence toward the asymptotic limits are different; and $m = o(n)$ if any $n_i \to \infty$.

Hall and Maiti (2006b) obtain that if $\hat{I}_\alpha$ is calibrated once (twice), the coverage accuracy improves to $O(m^{-2})$ $[O(m^{-3})]$. If the interval in (3.23) or (3.24) is calibrated once (twice), the coverage accuracy improves to $O(n^{-5/2})$ $[O(n^{-7/2})]$ when parameter dimension is fixed.

In summary, Hall and Maiti (2006b) cover a wide ranging independent data framework, with careful MSPE estimation and marginal coverage of prediction intervals as number of small areas increases; while we consider deeper linear mixed framework allowing for longitudinal dependence, and establish results as total data size increases on the joint variability of random effects and data, thus also obtaining area specificity.

**4. A simulation example.** In this section we compare the performance of our proposed parametric bootstrap with that of the traditional approaches, using a simulation study. We have carried out more extensive simulations which reflect the general pattern of performance reported here; the details are available from the authors.

TABLE 1
*Average coverage and average length of different intervals (nominal coverage = 0.95) in simulation pattern* (a)

| Group | Cox | FH | PR | PB–ET | PB–SL |
|---|---|---|---|---|---|
| G1 | 83.1 (3.12) | 90.4 (3.57) | 92.4 (3.82) | 96.1 (4.50) | 95.7 (4.42) |
| G2 | 85.4 (2.14) | 93.7 (2.50) | 98.0 (3.19) | 96.2 (2.83) | 95.9 (2.79) |
| G3 | 85.8 (2.02) | 93.9 (2.36) | 98.0 (3.08) | 96.0 (2.65) | 95.6 (2.61) |
| G4 | 86.1 (1.89) | 94.3 (2.19) | 98.2 (2.93) | 96.1 (2.43) | 95.7 (2.39) |
| G5 | 89.7 (1.12) | 95.2 (1.23) | 97.3 (1.87) | 95.7 (1.28) | 95.3 (1.26) |

For the sake of comparability with existing studies, we adopt part of the simulation framework of Datta, Rao and Smith (2005) for our study. We consider the Fay–Herriot model with $m = 15$ and $\mathbf{x}_i^T \beta = 0$, and consider five groups of small areas with three areas in each group. Within each group, the $D_i$'s remain the same. There are two different patterns for the $D_i$'s: (a) 0.2, 0.4, 0.5, 0.6, 4.0 [this is pattern (c) of Datta, Rao and Smith (2005)] and (b) 0.4, 0.8, 1.0, 1.2, 8.0. For pattern (a), we took $A = 1$ in order to make the results comparable to Datta, Rao and Smith (2005). For pattern (b), we took $A = 2$ in order to make the variances twice that of pattern (a), but preserve the $B_i = D_i/(A + D_i)$ ratios.

We obtain all the results based on 10,000 simulation runs. The Prasad–Rao method-of-moments, and the Fay–Herriot method of estimating the variance component $A$ are considered. Tables 1 and 2 report the simulated coverage probabilities and average lengths of several different prediction intervals (with nominal coverage 0.95) under patterns (a) and (b), respectively. We consider three prediction intervals of the type $EBLUP \pm 1.96\sqrt{mspe}$, where *mspe* is an estimator of the MSPE of EBLUP. The Cox interval, discussed in Section 2, is obtained by using Prasad–Rao method-of-moment estimator of $A$. The Prasad–Rao (PR) interval estimator is obtained using that estimator of $A$ along with the Prasad–Rao (1990) MSPE estimator. The Fay–Herriot (FH) interval estimator is obtained by using Fay–Herriot method of moments estimator of $A$ [see Fay and Herriot (1979)], and the MSPE estimator of EBLUP considered by Datta, Rao and Smith (2005). Along with these three, we report two different parametric bootstrap prediction intervals. In both the methods, we used Fay–Herriot method of estimating $A$, and the weighted least squares estimator of $\beta$. The first bootstrap interval is equal-tailed (PB–ET), and the second is the shortest length prediction interval (PB–SL). For both cases, we considered bootstrap sample of size 1000.

TABLE 2
*Average coverage and average length of different intervals (nominal coverage = 0.95) in simulation pattern* (b)

| Group | Cox | FH | PR | PB–ET | PB–SL |
|---|---|---|---|---|---|
| G1 | 85.5 (4.87) | 89.5 (5.18) | 89.3 (5.35) | 95.7 (6.55) | 95.4 (6.47) |
| G2 | 83.6 (2.68) | 86.0 (2.82) | 87.3 (2.93) | 95.2 (3.80) | 94.9 (3.75) |
| G3 | 83.4 (2.49) | 85.7 (2.60) | 86.8 (2.71) | 95.2 (3.53) | 94.9 (3.49) |
| G4 | 82.9 (2.27) | 85.0 (2.36) | 86.2 (2.46) | 95.0 (3.22) | 94.5 (3.18) |
| G5 | 83.0 (1.21) | 84.0 (1.23) | 84.8 (1.29) | 94.9 (1.72) | 94.6 (1.70) |

The figures in Table 1 are average coverage probabilities and average lengths for each prediction interval method for pattern (a), average being taken over all three small areas within each group. Table 2 reports similar results for pattern (b). It is clear that the results depend on the pattern of $D_i$'s. The Cox prediction interval method consistently undercover. For pattern (a), both parametric bootstrap prediction interval methods perform better than the Prasad–Rao and Fay–Herriot prediction intervals in terms of coverage errors. In this case, the Fay–Herriot method interval always undercovers, while the Prasad–Rao method interval switches from undercoverage to considerable over-coverage. The Prasad–Rao and the Fay–Herriot methods suffer from large undercoverage errors for pattern (b). In contrast, the performances of our parametric bootstrap methods remain stable over these two different patterns and always close to the target nominal level. Our minimum length parametric bootstrap method tends to provide shorter prediction intervals compared to the equal-tailed equivalents.

These performance patterns are repeated for other sample sizes, and other patterns of $D_i$ values in our simulations. It is generally seen that an increase in the variances results in the traditional intervals performing even more poorly. Thus, while both the Prasad–Rao and the Fay–Herriot MSPE estimators enjoy good theoretical properties, the resulting interval estimates suffer owing to the enforced symmetry and normality assumption.

## APPENDIX

PROOF OF THEOREM 3.1. We establish this result by obtaining an asymptotic expansion of $\mathcal{L}_n(q)$. An identical expansion holds for $\mathcal{L}_n^*(q)$, which leads to the result. In this proof, the letter capital $C$, with or without suffix, will be generic for constants.

For the projection on the column space of $\mathbf{X}$ we use the notation $\mathbf{P_x}$, thus
$$\mathbf{P_x} = \mathbf{X}(\mathbf{X}^T\mathbf{X})^{-1}\mathbf{X}^T.$$

Let $\phi(\cdot)$ ($\Phi(\cdot)$) be the standard Normal probability density (cumulative distribution) function. Let $\phi'$ and $\phi''$ denote the first and second derivative of $\phi(\cdot)$, thus for $x \in \mathbb{R}$ we have $\phi'(x) = -x\phi(x), \phi''(x) = (x^2 - 1)\phi(x)$. Define
$$Q(q, \mathbf{Y}) = \sigma_T^{-1}\{\hat{\mu}_T - \mu_T + q(\hat{\sigma}_T - \sigma_T)\}.$$

Then for any $q \in \mathbb{R}$, we have

$$\begin{aligned}
\mathcal{L}_n(q) &= \mathbb{P}(\hat{\sigma}_T^{-1}(T - \hat{\mu}_T) \leq q) \\
&= \mathbb{E}[\mathbb{P}(\sigma_T^{-1}(T - \mu_T) \leq q + \sigma_T^{-1}\{\hat{\mu}_T - \mu_T + q(\hat{\sigma}_T - \sigma_T)\}|\mathbf{Y})] \\
&= \mathbb{E}[\Phi(q + Q(q, \mathbf{Y}))] \\
&= \Phi(q) + \phi(q)\mathbb{E}Q(q, \mathbf{Y}) - 2^{-1}q\phi(q)\mathbb{E}Q^2(q, \mathbf{Y}) \\
&\quad + 2^{-1}\mathbb{E}\bigg\{\int_q^{q+Q}(q + Q - x)^2(x^2 - 1)\phi(x)\,dx\bigg\} \\
&= \Phi(q) + \phi(q)T_1(q) - 2^{-1}q\phi(q)T_2(q) + T_3(q).
\end{aligned}$$

Notice that for $x \in (q, q + Q)$, we have $0 \leq |q + Q - x| \leq |Q|$ and $(x^2 - 1) \times \phi(x) \leq 2\phi(\sqrt{3})$, we have

$$\begin{aligned}
\mathbb{E}\int_q^{q+Q}(q + Q - x)^2(x^2 - 1)\phi(x)\,dx \\
\leq \mathbb{E}\int_q^{q+Q}|(q + Q - x)^2||(x^2 - 1)\phi(x)|\,dx \\
\leq \mathbb{E}Q^2\int_q^{q+Q}2\phi(\sqrt{3})\,dx \leq C\mathbb{E}|Q|^3.
\end{aligned}$$

From the following calculations it will follow that $\mathbb{E}Q^8 = O(d^8 n^{-4})$, whereby $\sup_q T_3(q) = O(d^3 n^{-3/2})$.

We now simplify the expression for $Q(q, \mathbf{Y})$. Let $\hat{R} = R(\hat{\psi}), \hat{D} = D(\hat{\psi}), \hat{\Sigma} = \Sigma(\hat{\psi})$. Note that $\hat{R}\hat{\Sigma}^{-1} + \mathbf{Z}\hat{D}\mathbf{Z}^T\hat{\Sigma}^{-1} - R\Sigma^{-1} - \mathbf{Z}D\mathbf{Z}^T\Sigma^{-1} = 0$ almost surely. Hence we have

$$\begin{aligned}
\hat{\mu}_T - \mu_T &= c^T\hat{R}\hat{\Sigma}^{-1}\mathbf{X}\hat{\beta} + c^T\mathbf{Z}\hat{D}\mathbf{Z}^T\hat{\Sigma}^{-1}\mathbf{Y} - c^T R\Sigma^{-1}\mathbf{X}\beta - c^T\mathbf{Z}D\mathbf{Z}^T\Sigma^{-1}\mathbf{Y} \\
&= c^T[\mathbf{I} - \mathbf{Z}D\mathbf{Z}^T\Sigma^{-1}]\mathbf{P_x}(\mathbf{Zv} + \mathbf{e}) \\
&\quad + c^T(\mathbf{Z}\hat{D}\mathbf{Z}^T\hat{\Sigma}^{-1} - \mathbf{Z}D\mathbf{Z}^T\Sigma^{-1})(\mathbf{I} - \mathbf{P_x})(\mathbf{Zv} + \mathbf{e}) \quad \text{almost surely.}
\end{aligned}$$

In view of the above, let us write $Q(q, \mathbf{Y}) = Q_1 + Q_2(\mathbf{Y}) + Q_3(q, \mathbf{Y})$, where

$$\begin{aligned}
Q_1 &= \sigma_T^{-1}c^T[\mathbf{I} - \mathbf{Z}D\mathbf{Z}^T\Sigma^{-1}]\mathbf{P_x}(\mathbf{Zv} + \mathbf{e}), \\
Q_2(\mathbf{Y}) &= \sigma_T^{-1}c^T(\mathbf{Z}\hat{D}\mathbf{Z}^T\hat{\Sigma}^{-1} - \mathbf{Z}D\mathbf{Z}^T\Sigma^{-1})(\mathbf{I} - \mathbf{P_x})(\mathbf{Zv} + \mathbf{e}), \\
Q_3(q, \mathbf{Y}) &= q\sigma_T^{-1}(\hat{\sigma}_T - \sigma_T).
\end{aligned}$$

First, using the decomposition

$$Q_1 = \sigma_T^{-1}c^T\mathbf{P_x}(\mathbf{Zv} + \mathbf{e}) - \sigma_T^{-1}c^T\mathbf{Z}D\mathbf{Z}^T\Sigma^{-1}\mathbf{P_x}(\mathbf{Zv} + \mathbf{e}) = Q_{11} - Q_{12}.$$

Using Assumption 1, in particular, (3.4)–(3.7), with some amount of algebra we can conclude that $\mathbb{E}Q_1 = 0$, $\mathbb{E}Q_{11}^2 = O(n^{-1})$, $\mathbb{E}Q_{11}^8 = O(p^4 n^{-4})$, $\mathbb{E}Q_{12}^2 = O(n^{-1})$, $\mathbb{E}Q_{12}^8 = O(p^4 n^{-4})$.

We now analyze $Q_2(\mathbf{Y})$ and $Q_3(q, \mathbf{Y})$. These are considerably more complicated than $Q_1$. We initially break down these two quantities in terms of more tractable variables $W_1, \ldots, W_{11}$ and remainder terms. The variables $W_1, \ldots, W_{11}$ depend on $\mathbf{Z}\hat{D}\mathbf{Z}^T$ and $\hat{\Sigma}^{-1}$ and their population equivalents. We need to compute the first, second, fourth, eighth and sixteenth moment of the $W_i$'s and show that the remainder terms are negligible.

Toward that goal, our next step is to expand $\mathbf{Z}\hat{D}\mathbf{Z}^T$ in equation (A.1) and $\hat{\Sigma}^{-1}$ in equation (A.2) in terms of simpler matrices. Then we obtain asymptotic expansions of the matrix entries, whereby at last we have sufficient ingredients for the moment computations of $W_1, \ldots, W_{11}$ and the remainder terms. We skip the details of the moment calculation algebra, of which there are several hundreds to compute. However, our assumptions are sufficient to establish the end result that $\mathbb{E}Q_2(\mathbf{Y})$, $\mathbb{E}Q_2^2(\mathbf{Y})$, $\mathbb{E}Q_3(q, \mathbf{Y})$ and $\mathbb{E}Q_3^2(q, \mathbf{Y})$ are all $O(d^2/n)$. The expansion of $Q_2(\mathbf{Y})$ is as follows:

$$\begin{aligned}
Q_2(\mathbf{Y}) &= \sigma_T^{-1} c^T (\mathbf{Z}\hat{D}\mathbf{Z}^T \hat{\Sigma}^{-1} - \mathbf{Z}D\mathbf{Z}^T \Sigma^{-1})(\mathbf{I} - \mathbf{P_x})(\mathbf{Z}\mathbf{v} + \mathbf{e}) \\
&= \sigma_T^{-1} c^T [(\mathbf{Z}\hat{D}\mathbf{Z}^T - \mathbf{Z}D\mathbf{Z}^T)\Sigma^{-1} \\
&\quad + \mathbf{Z}D\mathbf{Z}^T (\hat{\Sigma}^{-1} - \Sigma^{-1}) \\
&\quad + (\mathbf{Z}\hat{D}\mathbf{Z}^T - \mathbf{Z}D\mathbf{Z}^T)(\hat{\Sigma}^{-1} - \Sigma^{-1})](\mathbf{I} - \mathbf{P_x})(\mathbf{Z}\mathbf{v} + \mathbf{e}) \\
&= W_1 + W_2 + W_3.
\end{aligned}$$

Now define $W = \sigma_T^{-2}(\hat{\sigma}_T^2 - \sigma_T^2)$. We will simplify $Q_3(q, \mathbf{Y})$ in terms of $W$. However, we first need to simplify $W$. For this, we have

$$\begin{aligned}
W &= \sigma_T^{-2}(\hat{\sigma}_T^2 - \sigma_T^2) \\
&= \sigma_T^{-2} c^T [\{\mathbf{Z}\hat{D}\mathbf{Z}^T - \mathbf{Z}D\mathbf{Z}^T\} \\
&\quad - \{\mathbf{Z}\hat{D}\mathbf{Z}^T \hat{\Sigma}^{-1} \mathbf{Z}\hat{D}\mathbf{Z}^T - \mathbf{Z}D\mathbf{Z}^T \Sigma^{-1} \mathbf{Z}D\mathbf{Z}^T\}] c \\
&= \sigma_T^{-2} c^T [\{\mathbf{Z}\hat{D}\mathbf{Z}^T - \mathbf{Z}D\mathbf{Z}^T\} - \mathbf{Z}D\mathbf{Z}^T (\hat{\Sigma}^{-1} - \Sigma^{-1})\mathbf{Z}D\mathbf{Z}^T \\
&\quad - \mathbf{Z}D\mathbf{Z}^T \Sigma^{-1} \{\mathbf{Z}\hat{D}\mathbf{Z}^T - \mathbf{Z}D\mathbf{Z}^T\} - \{\mathbf{Z}\hat{D}\mathbf{Z}^T - \mathbf{Z}D\mathbf{Z}^T\}\Sigma^{-1}\mathbf{Z}D\mathbf{Z}^T \\
&\quad - \{\mathbf{Z}\hat{D}\mathbf{Z}^T - \mathbf{Z}D\mathbf{Z}^T\}\Sigma^{-1}\{\mathbf{Z}\hat{D}\mathbf{Z}^T - \mathbf{Z}D\mathbf{Z}^T\} \\
&\quad - \{\mathbf{Z}\hat{D}\mathbf{Z}^T - \mathbf{Z}D\mathbf{Z}^T\}(\hat{\Sigma}^{-1} - \Sigma^{-1})\mathbf{Z}D\mathbf{Z}^T \\
&\quad - \mathbf{Z}D\mathbf{Z}^T (\hat{\Sigma}^{-1} - \Sigma^{-1})\{\mathbf{Z}\hat{D}\mathbf{Z}^T - \mathbf{Z}D\mathbf{Z}^T\} \\
&\quad - \{\mathbf{Z}\hat{D}\mathbf{Z}^T - \mathbf{Z}D\mathbf{Z}^T\}(\hat{\Sigma}^{-1} - \Sigma^{-1})\{\mathbf{Z}\hat{D}\mathbf{Z}^T - \mathbf{Z}D\mathbf{Z}^T\}] c
\end{aligned}$$

$$= \sum_{i=4}^{11} W_i.$$

Let us now simplify $Q_3(q, \mathbf{Y})$:

$$Q_3(q, \mathbf{Y}) = q\sigma_T^{-1}(\hat{\sigma}_T - \sigma_T) = q(\sigma_T^{-1}\hat{\sigma}_T - 1)$$
$$= q[W/2 - W^2/8 + r_n].$$

At this stage, we use a result from Rao (1965), page 41, result (1c.3.10): If $A$ is positive definite and $B$ nonnegative definite $n \times n$ symmetric matrices, then, we can write $A = \sum_{i=1}^n A_i A_i^T$, $B = \sum_{i=1}^n b_i A_i A_i^T$ where $A_i$'s are vectors, $b_i$'s are real constants. Moreover, since $A$ is positive definite, the $A_i$'s form a basis (but not necessarily an orthogonal basis) of $\mathbb{R}^n$. Another way of writing the same thing is $A = A_1 A_1^T$, $B = A_1 \Lambda_B A_1^T$ where the columns of $A_1$ are the $A_i$'s, and $\Lambda_B$ is a diagonal matrix with entries $b_i$. Note that $A_1$ is nonsingular.

We use this result twice. First, we take $R$ as $A$ and $\hat{R}$ as $B$, and then we take $D$ as $A$ and $\hat{D}$ as $B$. Thus we have

$$R = R_1(\psi)R_1^T(\psi), \qquad \hat{R} = R_1(\psi)\Lambda_R(\hat{\psi})R_1^T(\psi),$$
$$D = D_1(\psi)D_1^T(\psi), \qquad \hat{D} = D_1(\psi)\Lambda_D(\hat{\psi})D_1^T(\psi).$$

Here the nonsingular matrices $R_1$ and $D_1$ depend on the unknown parameter $\psi$, while $\Lambda_R$ and $\Lambda_D$ are diagonal matrices depending on the estimator $\hat{\psi}$.

Based on the above, we have

$$\Sigma = R_1[\mathbf{I} + R_1^{-1}\mathbf{Z}D_1 D_1^T \mathbf{Z}^T R_1^{-T}]R_1^T = R_1[\mathbf{I} + AA^T]R_1^T,$$

where $A = R_1^{-1}\mathbf{Z}D_1$. We define $B_0 = [\mathbf{I} + AA^T]^{-1/2}$, the symmetric square root. Hence, $\Sigma^{-1} = R_1^{-T} B_0^2 R_1^{-1}$. We also have

$$\hat{\Sigma} = R_1[\Lambda_R + R_1^{-1}\mathbf{Z}D_1\Lambda_D D_1^T \mathbf{Z}^T R_1^{-T}]R_1^T = R_1[\Lambda_R + A\Lambda_D A^T]R_1^T.$$

Let us write $\Lambda_R = \mathbf{I} + (k/n)^{1/2}U_R$, $\Lambda_D = \mathbf{I} + (k/n)^{1/2}U_D$. Our next step is to write $\mathbf{Z}\hat{D}\mathbf{Z}^T$ and $\hat{\Sigma}^{-1}$ using $U_R$ and $U_D$. Thus

$$\mathbf{Z}\hat{D}\mathbf{Z}^T = R_1 A \Lambda_D A^T R_1^T$$
$$= \mathbf{Z}D\mathbf{Z}^T + (k/n)^{1/2} R_1 A U_D A^T R_1^T,$$

(A.1)

$$\hat{\Sigma}^{-1} = R_1^{-T}[\Lambda_R + A\Lambda_D A^T]^{-1} R_1^{-1}$$
$$= R_1^{-T} B_0 [\mathbf{I} + (k/n)^{1/2} U]^{-1} B_0 R_1^{-1},$$

where $U = B_0(U_R + AU_D A^T)B_0$. We further simplify $\hat{\Sigma}^{-1}$ by writing

$$[\mathbf{I} + (k/n)^{1/2}U]^{-1}$$
$$= \mathbf{I} - (k/n)^{1/2}U + (k/n)U^2 - (k/n)^{3/2}[\mathbf{I} + (k/n)^{1/2}U]^{-1}U^3$$
$$= \mathbf{I} - (k/n)^{1/2}U + (k/n)U^2 - (k/n)^{3/2}U_R.$$

Hence,

$$\hat{\Sigma}^{-1} = R_1^{-T} B_0 [\mathbf{I} + (k/n)^{1/2} U]^{-1} B_0 R_1^{-1},$$

(A.2)
$$= \Sigma^{-1} - (k/n)^{1/2} R_1^{-T} B_0 U B_0 R_1^{-1} + (k/n) R_1^{-T} B_0 U^2 B_0 R_1^{-1}$$
$$- (k/n)^{3/2} R_1^{-T} B_0 U_R B_0 R_1^{-1}.$$

Equations (A.1) and (A.2) will be heavily used in the analysis below.

We now turn our attention to $U_R$ and $U_D$. Recall that $S = k^{1/2} n^{-1/2} (\hat{\psi} - \psi)$. Suppose $\lambda_{Di}$ is the $i$th element of either $\Lambda_D$. We have $\lambda_{Di}(\psi) = 1$. Thus, we have

$$\hat{\lambda}_{Di} = \lambda_{Di}(\hat{\psi})$$
$$= 1 + (k/n)^{1/2} S^T \lambda'_{Di} + 2^{-1} (k/n) S^T \lambda''_{Di} S$$
$$+ 6^{-1} (k/n)^{3/2} \sum_{j_1, j_2, j_3} \Delta(j_1, j_2, \ldots, j_3; \lambda_{Di}(\psi^*)) S_{j_1} S_{j_2} S_{j_3},$$

where $\psi^*$ is a point between $\psi$ and $\hat{\psi}$. Hence, we have

$$(k/n)^{1/2} U_{Di} = (\lambda_{Di}(\hat{\psi}) - 1)$$
$$= (k/n)^{1/2} S^T \lambda'_{Di} + 2^{-1} (k/n) S^T \lambda''_{Di} S$$
$$+ 6^{-1} (k/n)^{3/2} \sum_{j_1, j_2, j_3} \Delta(j_1, j_2, \ldots, j_3; \lambda_{Di}(\psi^*)) S_{j_1} S_{j_2} S_{j_3}$$
$$= (k/n)^{1/2} U_{Di1} + (k/n) U_{Di2} + (k/n)^{3/2} U_{Di3}.$$

A similar analysis holds for $U_{Ri}$.

It now remains to calculate the first, second and eighth moments of $W_1, \ldots, W_{11}$, and establish that $\mathbb{E} W_i^8 = O(k^{16} n^{-8})$, $\mathbb{E} W_i^2 = O(k^2/n)$, $\mathbb{E} W_i = O(k/n)$ for $i = 1, \ldots, 11$. Some of these moments turn out to be of even smaller order and thus contribute negligibly. Also, certain remainder terms have to proved negligible. The totality of these computations involve a few hundred algebraic manipulations, and is reasonably routine. We sketch part of the computation for one of the components of $W_1$ as an example of the technique used. The rest of the computations are omitted.

Using (A.1) we obtain that

$$W_1 = \sigma_T^{-1} c^T (\mathbf{Z} \hat{D} \mathbf{Z}^T - \mathbf{Z} D \mathbf{Z}^T) \Sigma^{-1} (\mathbf{I} - \mathbf{P_x})(\mathbf{Z} \mathbf{v} + \mathbf{e})$$
$$= (k/n)^{1/2} \sigma_T^{-1} c^T \mathbf{Z} D_1$$
$$\times [U_{D1} + (k/n)^{1/2} U_{D2} + (k/n) U_{D3}] D_1 \mathbf{Z}^T \Sigma^{-1} (\mathbf{I} - \mathbf{P_x})(\mathbf{Z} \mathbf{v} + \mathbf{e})$$
$$= W_{11} + W_{12} + W_{13}.$$

Here $U_{Dj}$, $j = 1, 2, 3$ is the diagonal matrix whose $i$th diagonal entry is $U_{Dij}$ computed earlier. We will sketch the computations only for

$$W_{11} = (k/n)^{1/2} \sigma_T^{-1} c^T \mathbf{Z} D_1 U_{D1} D_1 \mathbf{Z}^T \Sigma^{-1} (\mathbf{I} - \mathbf{P_x})(\mathbf{Zv} + \mathbf{e}).$$

Recall that $U_{Di1} = S^T \lambda'_{Di}$. Let the $j$th element of the $\mathbb{R}^k$ dimensional vector $\lambda_{Di}$ be $\lambda'_{Dij}$, $j = 1, \ldots, k$, $i = 1, \ldots, q$. Define the $(k \times q)$ matrix $\Lambda'_D$ whose $(a, b)$th element is $\lambda'_{Dba}$. Also define the diagonal matrix $E_3$ whose $i$th diagonal entry is the $i$th element of the vector $D_1 \mathbf{Z}^T c$. Then note that $W_{11} = (k/n)^{1/2} \sigma_T^{-1} S^T E_5 (\mathbf{Zv} + \mathbf{e})$, where $U \sim N_n(0, \mathbf{I}_n)$, $S$ depends on $U$, and $E_5 = \Lambda'_D E_3 D_1 \mathbf{Z}^T \Sigma^{-1} (\mathbf{I} - \mathbf{P_x})$. The appropriate moment properties of $W_{11}$ now follow by applying (3.20), (3.21) and (3.12).

The above sketch of calculations for $W_{11}$ may be repeated with variations for the other terms as well to establish that

$$\mathcal{L}_n(q) = \Phi(q) + k^2 n^{-1} \gamma(q, \beta, \psi) + O(k^3 n^{-3/2}),$$

for a $O(1)$ smooth quantity $\gamma(\cdot, \cdot, \cdot)$.

A similar representation holds for $\mathcal{L}_n^*(q)$ with $\hat{\beta}$ and $\hat{\psi}$ in place of $\beta$ and $\psi$. This establishes the result. $\square$

S. Chatterjee  
School of Statistics  
University of Minnesota  
224 Church Street SE  
Minneapolis, Minnesota 55455  
USA  
E-mail: chatterjee@stat.umn.edu  
URL: http://www.stat.umn.edu/~chatterjee/

P. Lahiri  
H. Li  
Department of Mathematics  
University of Maryland, College Park  
College Park, Maryland 20742  
USA  
E-mail: plahiri@survey.umd.edu  
        huilin@math.umd.edu  
URL: http://www.jpsm.umd.edu